\documentclass[11pt]{article} 

\usepackage{amssymb}

\usepackage{babel}
\usepackage[latin1]{inputenc}
\usepackage[OT1]{fontenc} 

\usepackage{epic, eepic}

\newtheorem{thm}{Theorem}
\newtheorem{theorem}[thm]{Theorem}
\newtheorem{definition}[thm]{Definition}
\newtheorem{lemma}[thm]{Lemma}

\newtheorem{corollary}[thm]{Corollary}

\newtheorem{remark}[thm]{Remark}

\newcommand{\pf}{{\bf Proof: \/}}

\newcommand\qed{\hfill{\nobreak\quad\raise -2pt\hbox{\vrule\vbox to
      8pt{\hrule width 8pt \vfill\hrule}\vrule}\par\vspace{2ex}}}

\def\emm#1,{{\em#1}\/}
\newcommand\ki{\raise.485ex\hbox{$\chi$}}
\def\ch#1,#2,{{#1\choose #2}}
\newcommand\NN{\mathbb N}
\newcommand\cl[1]{{\mc #1}}
\newcommand\mc{\mathcal}

\newcommand\sm{\ensuremath{\setminus}}
\newcommand\dg{\ensuremath{\Delta_G}}
\newcommand\kig{\ensuremath{\ki_G}}
\newcommand\tg{\ensuremath{T_G}}
\newcommand\bdh{\ensuremath{\tilde{B}_{d}}}

\def\disp{\displaystyle}
\newcommand\cls{{\mc S}}
\newcommand\st{\; | \;}
\newcommand\sseq{\subseteq}
\newcommand\ssneq{\subsetneq}

\newcommand\lra{\longrightarrow}
\newcommand\cd{\cdot}
\newcommand{\pref}[1]{{(\protect\ref{#1})}}

\begin{document}
\thispagestyle{empty}

\begin{center} {\Large \bf The Coloring Ideal and Coloring Complex of a Graph}\\

  \vspace{12pt}
  {\large Einar Steingr\'{\i}msson\\
    {Matematik CTH \& GU}\\ {412 96 G\"oteborg}}\\
  {SWEDEN}\\
  {\normalsize\tt einar@math.chalmers.se}\\
\end{center}

\begin{abstract}
  Let $G$ be a simple graph on $d$ vertices.  We define a monomial
  ideal $K$ in the Stanley-Reisner ring $A$ of the order complex of
  the Boolean algebra on $d$ atoms.  The monomials in $K$ are in one-to-one
  correspondence with the proper colorings of $G$.  In particular, the
  Hilbert polynomial of $K$ equals the chromatic polynomial of $G$.
  
  The ideal $K$ is generated by square-free monomials, so $A/K$ is the
  Stanley-Reisner ring of a simplicial complex $C$.  The $h$-vector of
  $C$ is a certain transformation of the \emm tail, $T(n)=
  n^d-\ki(n)$ of the chromatic polynomial $\ki$ of $G$.
  The combinatorial structure of the complex $C$ is described
  explicitly and it is shown that the Euler characteristic of $C$
  equals the number of acyclic orientations of $G$.
\end{abstract}

\section{Introduction and preliminaries}

Let $G$ be a simple graph on $d$ vertices.  In this paper we construct
a monomial ideal $K$ in the face ring (Stanley-Reisner ring) $A$ of
the order complex of the Boolean algebra on $d$ atoms, which is
isomorphic to a cone over the first barycentric subdivision of a
$(d-1)$-simplex.  The monomials in $K$ are in one-to-one
correspondence with the proper colorings of $G$.  The quotient of $A$
by $K$ is the face ring of a simplicial complex whose structure can be
described explicitly.

The construction of the ideal $K$ is based on a definition of Chung
and Graham \cite{coverpoly}, whose purpose was to give a combinatorial
interpretation to the coefficients of the chromatic polynomial
$\ki(n)$ of $G$ when written in the basis $\{\ch
n+k,d,\}_{k=0,\ldots,d}$.  It was shown by Chow \cite{chow-quasi} that
this result can also be derived from a theorem of Stanley's concerning
his \emm chromatic symmetric function, \cite{stan-symm}.  However, our
construction does not seem to be much related to Stanley's chromatic
symmetric function.  In fact, our complexes are isomorphic for the two
non-isomorphic graphs on $n\ge4$ vertices and two edges, whereas
Stanley's functions for these graphs are not equal.  On the
other hand, our complex \emm does, distinguish the two non-isomorphic
(but chromatically equivalent) graphs on five vertices that Stanley's
function does not distinguish.

As it happens, invariants related to colorings of a graph $G$ often
have connections to the acyclic orientations of $G$.  In our case, we
show that the Euler characteristic of the coloring complex equals the
number of acyclic orientations of $G$.

It should also be mentioned that Brenti asked \cite{brenti-hilbert}
whether there exists, for an arbitrary graph $G$, a standard graded
algebra whose Hilbert polynomial equals the chromatic polynomial of
$G$.  This was answered in the affirmative by Almkvist
\cite{almkvist}, but his proof is non-constructive.  The structure of
such a graded algebra, however, will not necessarily be closely
related to the colorings of $G$, since its monomials of degree less
than the chromatic number of $G$ cannot correspond to colorings of
$G$.

Throughout this paper, unless otherwise specified, $G$ is a graph on
$d$ vertices labeled by the elements of $[d] = \{1,2,\ldots,d\}$, with
no loops and no multiple edges.  Frequently, but not always, we
suppress $G$ from the notation for simplicity.

A \emm path of length $k$, in $G$ is a sequence
$v_0,v_1,v_2,\ldots,v_k$ of vertices of $G$ such that there is an edge
between $v_{i-1}$ and $v_{i}$, for each $i\in[k]$.

A \emm stable, set in $G$ is a set of vertices with no edge between
any pair.

Let $V$ be the set of vertices of $G$.  A \emm coloring of $G$, is a
map $\phi:V\lra\NN$ with $\phi(x)\ne\phi(y)$ if $x$ and $y$ are
adjacent in $G$, that is, if there is an edge between $x$ and $y$.
Thus, we treat the natural numbers as colors and when referring to the
ordering of colors, we mean the usual ordering on $\NN$.  (Observe
that we omit the word ``proper'' from the definition of coloring,
since we only consider proper colorings).

A \emm coloring of $G$ with $n$ colors,, or \emm $n$-coloring,, is a
map $\phi:V\lra[n]$ with $\phi(x)\ne\phi(y)$ if $x$ and $y$ are
adjacent.  (Observe that $\phi$ need not be surjective.)

\begin{definition}\rm%
  Let $S_1,S_2,\ldots,S_m$ be an ordered partition of the vertices of
  $G$.  For $v\in G$ let $\ell(v)$ be the length of the longest path
  $v_{i_1},v_{i_2},\ldots,v_{i_p}=v$ (ending in $v$) in $G$ such that
  $v_{i_j}\in S_{i_j}$ for each $j$ and $i_1<i_2<\cdots<i_p$.
  
  If $\pi= a_1 a_2\cdots a_d$ is a permutation in the symmetric group
  $\cl S_d$ we let $\pi$ induce the ordered partition of the vertices
  in $G$ obtained by letting $a_i$ constitute the $i$-th block
  (singleton) in the ordered partition.  In accordance with the
  definition of $\ell(v)$ subject to an ordered partition, we then let
  $\ell(k)$ be the length of the longest path $a_{i_1},
  a_{i_2},\ldots, a_{k}$ (ending in $a_k$) in $G$ such that
  $i_1<i_2<\cdots<k$.
\end{definition} 

The following definition is a variation of a definition of Chung and
Graham in \cite[\S 5]{coverpoly}.

\begin{definition}\rm\label{cut-def}%
  The integer $k\in[0,\ldots,d-1]$ is a \emm cut in $\pi$ (with respect to
  $G$), if
\begin{enumerate}%
\item%
$k =0$, or
\item%
$\ell(k) < \ell(k+1)$, or
\item%
$\ell(k) = \ell(k+1)$ and $a_k<a_{k+1}$.
\end{enumerate}
\end{definition} 

\begin{definition}\rm\label{seq-def}%
  Let $\pi=a_1 a_2\cdots a_d$ be a permutation with cuts
  $i_1=0,i_2,\ldots,i_k$. The \emm $G$-sequence, of $\pi$ is the
  sequence of sets $S_{1},S_{2},\ldots,S_{k}$ where $S_j=\{a_{i_j
    +1},a_{i_j +2},\ldots, a_{i_{j +1}}\}$, for $j<k$, and
  $S_k=\{a_{i_k +1},a_{i_k +2},\ldots, a_d\}$.  The \emm short
  $G$-sequence of $\pi$, is $S_{1},S_{2},\ldots,S_{k-1}$.
\end{definition}

As an example, let $G$ be the graph
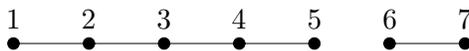
\setlength{\unitlength}{1mm}
%
\begin{center}
\begin{picture}(50,10)(0,-3)
\newcommand\pt{\circle*{1.5}}

\path(0,0)(40,0)
\path(50,0)(60,0)

\put(0,0){\pt}
\put(10,0){\pt}
\put(20,0){\pt}
\put(30,0){\pt}
\put(40,0){\pt}
\put(50,0){\pt}
\put(60,0){\pt}

\put(0,0){\makebox(0,6.5){1}}
\put(10,0){\makebox(0,6.5){2}}
\put(20,0){\makebox(0,6.5){3}}
\put(30,0){\makebox(0,6.5){4}}
\put(40,0){\makebox(0,6.5){5}}
\put(50,0){\makebox(0,6.5){6}}
\put(60,0){\makebox(0,6.5){7}}
\end{picture}
\label{fig1}
\end{center}
%
and let $\pi$ be the permutation 5236417.  Then the path lengths
$\ell(k)$ associated to $\pi$ are given by
$$
\begin{array}{ccccccc}
5 & 2 & 3 & 6 & 4 & 1 & 7\cr
0 & 0 & 1 & 0 & 2 & 1 & 1
\end{array}
$$
so $\pi$ has cuts 0, 2, 4 and 6 and thus  $G$-sequence $\{2,5\}, ~
\{3,6\}, ~ \{1,4\}, ~ \{7\}$.

\begin{lemma}\label{stable}%
  Let $S_1,S_2,\ldots,S_k$ be the  $G$-sequence of a permutation
  $\pi$.  If $a_m, a_{m+1}\in S_i$ then either $\ell(m) > \ell(m+1)$ or
  else $\ell(m) = \ell(m+1)$ and $a_m>a_{m+1}$.  Moreover, each $S_i$
  is a stable set in $G$.
\end{lemma}
\pf{%
  The first part of the lemma follows directly from Definitions
  \ref{cut-def} and \ref{seq-def}.  Thus, if $a_m$ and $a_k$ both
  belong to $S_i$, with $m<k$, then $\ell(m)\ge\ell(k)$.  But if $a_m$
  and $a_k$ are adjacent in $G$ then $\ell(m) < \ell(k)$, since a
  path ending in $a_m$ can always be extended to end in $a_k$.  That
  is a contradiction, so no two elements in $S_i$ are adjacent in $G$.
  }
\qed{}

The following theorem is stated (in a different, but equivalent, form)
in \cite{coverpoly}, where it is claimed that it follows from the work
of Brenti in \cite{brenti-expansions}.  Also, in \cite{chow-digraph}
it is shown that the theorem follows from certain properties of the
chromatic symmetric function of Stanley \cite{stan-symm}.  We give a
proof that is different from both of these, but which better suits our
purposes.  First, a definition.

\begin{definition}\rm%
  Let $P(n)$ be a polynomial of degree $d$.  The \emm $W$-transform of $P$, is the
  polynomial $W$ defined by
$$
\sum_{n\ge 0}{P(n)t^n}= \frac{W(t)}{(1-t)^{d+1}}.
$$
\end{definition} 
That $W$ is a polynomial is easily shown, as is the fact that its
degree is at most equal to the degree of $P$.
\begin{theorem}\label{w-theorem}%
  Let $W_G(t)$ be the $W$-transform of $\ki_G$, that is,
$$
\sum_{n\ge 0}{\kig(n)t^n}= \frac{W_G(t)}{(1-t)^{d+1}}.
$$
Then we have
$$
W_G(t) = \sum_{\pi\in\cl S_d} {t^{c(\pi)}},
$$
where $c(\pi)$ is the number of cuts in $\pi$.
In particular, $\disp\sum_{\pi\in\cl S_d} {t^{c(\pi)}}$ is independent of
the labeling of $G$.
\end{theorem}
\pf{%
  Define $w_i$, for $0\le i\le d$, by setting $W_G(t) =
  w_0+w_1t+\cdots+w_dt^d$.  Then the claim is equivalent to saying
  that
$$
\kig(n) = \sum_{k=0}^{d}  {\ch n+k,d, w_{d-k} },
$$
or, in other words, that there are, for each permutation in $\cl
S_d$ with $(d-k)$ cuts, exactly $\ch n+k,d,$ ways to color $G$ with
$n$ colors.  Let $\pi$ be such a permutation.  By the latter part of
Lemma~\ref{stable}, there is an obvious way to associate a permutation
$\pi$ with $k$ cuts to a coloring using $k$ (out of $n$ available)
colors.  Namely, color the vertices corresponding to the letters
between the $(i-1)$-st and the $i$-th cut with color number $i$, and
the letters after the last cut with color $k$.

Since $\pi$ has $(d-k)$ cuts, there are exactly $k$ places between
adjacent letters in $\pi$ that do not correspond to cuts.  Let us
pick, among $n$ colors and $k$ non-cuts, exactly $d$ items,
say $i$ non-cuts and $(d-i)$ colors. Then we have chosen which $(d-i)$
colors to use, and which $i$ non-cuts to retain.  We now introduce
extra cuts at the remaining $(k-i)$ non-cuts, which means we have a
total of $(d-k) +(k-i)=(d-i)$ cuts, and thereby $(d-i)$ stable sets,
which get colored by the $(d-i)$ colors chosen, in the order
prescribed by $\pi$.
 
 Thus, we have associated $\ch n+k,d,$ colorings to each permutation
 $\pi\in\cl S_d$ with $(d-k)$ cuts.  
 
 For the converse, we need to show that each coloring of $G$ using
 some of $n$ available colors arises from a unique permutation,
 together with a choice of extra cuts, as described above.  Given such
 a coloring, partition the vertices of $G$ into blocks, each
 consisting of all vertices with like color, and order the blocks
 increasingly by color.  Within each block, order the vertices so that
 their corresponding maximal path lengths $\ell(m)$ are (weakly)
 decreasing, and so that the vertices are decreasing when two vertices
 have the same maximal path length associated to them.  This can
 always be done, because the path lengths of two vertices in the same
 block depend only on the vertices in preceding blocks.  Thus, writing
 the vertices in the order described we get a unique permutation $\pi$
 in $\cls_d$.
 
 By the construction of $\pi$, all of its cuts occur between blocks of
 the ordered partition $P$ from which $\pi$ was constructed.  Thus,
 the coloring from which $P$ was constructed arises from $\pi$
 together with the extra cuts (separating blocks in $P$) that are not
 cuts in $\pi$.  }\qed{}

\section{The coloring ideal}

The field $k$ in the following definition can be taken to be the
complex numbers.  Also, all rings are taken to be commutative.  Recall
that the short $G$-sequence of a permutation $\pi$ is the $G$-sequence
of $\pi$ take away the last set in the sequence.

For undefined terminology and background in what follows, see
\cite{CCA}.

Let $A=k\left[x_S \st S\sseq[d]\right]$, that is, $A$ is the
polynomial ring whose indeterminates correspond to all subsets of
$[d]$.  Throughout, $R$ will denote the \emm face ring, (or \emm
Stanley-Reisner ring,) of the order complex of the Boolean algebra on
$d$ atoms.  This ring is the quotient $A/I$, where $I =\{x_S x_T\st
S\nsubseteq T \mbox{ and } T\nsubseteq S\}$.  Thus, the monomials of
$R$ correspond precisely to those monomials $M=x_{S_1}x_{S_2} \cdots
x_{S_k}\in A$ for which $\emptyset\sseq S_1 \sseq S_2\sseq \cdots\sseq
S_k\sseq [d]$ (for some (unique) rearrangement of the indices).

\begin{definition}\rm%
  Let $M=x_{S_1} x_{S_2} \cdots x_{S_k}$ be a monomial in $A$ such
  that
$$
\emptyset\ssneq S_1 \ssneq S_2\ssneq \cdots\ssneq S_k \ssneq [d],
$$
where $\ssneq$ denotes strict inclusion.
Then $M$ is a \emm basic coloring monomial for $G$, if there is a
permutation $\pi\in\cl S_d$ with short $G$-sequence $S_1,S_2\sm
S_1,\ldots,S_{k}\sm S_{k-1}$.
  
  A nonzero monomial $M=x_{S_1} x_{S_2} \cdots x_{S_k}\in A$ is a \emm
  coloring monomial for $G$, if $M$ is divisible by a basic coloring
  monomial for $G$.
\end{definition} 

\begin{definition}\rm%
The \emm coloring ideal of $G$, is the ideal $K_G$ in $R$ generated by
all (equivalently by the basic) coloring monomials for G.
\end{definition} 

Note that a coloring monomial may not be divisible by $x_{[d]}$, since
no basic coloring monomial is divisible by $x_{[d]}$, as it is
constructed from a \emm short,\/ $G$-sequence.  That we choose to
define the coloring ideal in this way is due to a technicality which
will be explained in Remark \ref{short-remark}.

We shall now show that the monomials of degree $n$ in $K_G$ are in
one-to-one correspondence with the colorings of $G$ with $n+1$ colors.

It follows from the definitions of $R$ and $K_G$ that any monomial $M$
in $K_G$ can be written as 
$$
x_{S_1}^{e_1} x_{S_2}^{e_2} \cdots x_{S_k}^{e_k}
$$
such that $S_1 \ssneq S_2\ssneq \cdots\ssneq S_k \sseq [d]$ and
such that $S_i\sm S_{i-1}$ is a stable set in $G$ for each $i$.  Such
a monomial gives rise to a unique coloring of $G$ with $n$ colors
where $n$ is the degree of $M$.  Namely, if $S_1=\emptyset$, the
colors $1,2,\ldots,e_1$ are not used.  Otherwise, the vertices in
$S_1$ get color 1.  The vertices in $S_2\sm S_1$ get color $e_1+1$
and, in general, the vertices in $S_i\sm S_{i-1}$ get color
$e_1+e_2+\cdots+e_{i-1}+1$.  The vertices in $[d]\sm S_k$ get color
$\sum_i{e_i}+1$.  If $S_k=[d]$ and $e_k>0$ then the last $e_k$ colors
are not used (recall that in an $n$-coloring not all $n$ colors have
to be used).

As an example, suppose 
$$
M=x_{\emptyset}^2\cd x_{25}^3\cd x_{235}^2
$$
is a coloring monomial for $G$ (where we write 25 for the set
$\{2,5\}$ etc.).  In the corresponding coloring of $G$ with $8$
colors, the vertices 2 and 5 get color 3, the vertex 3 gets color 6
and all the remaining vertices (however many they are) get color $8$.
Multiplying $M$ by $x_{[d]}^e$ corresponds to regarding the coloring
in question as a coloring with $8+e$ colors.  Clearly, two different
monomials yield different colorings.  If they are identical except for
the exponent to $x_{[d]}$ they correspond to colorings with a
different number of colors, although each vertex gets the same color
in each of the two colorings.

Conversely, suppose we have a coloring of $G$ with $n$ colors (of
which not all have to be used).  Then each color is associated to a
stable set.  Order these sets increasingly by color and let $S_i$ be
the union of the first $i$ of them.  Then we can construct a
corresponding coloring monomial as we now explain by an example.
Suppose we have a $9$-coloring of $G$ where vertices 3 and 6 have
color 4, vertex 7 has color 6 and vertices 1, 2, 4 and 5 have color 7.
Then we have the sequence of sets $\{3,6\}$, $\{3,6,7\}$ and
$\{1,2,3,4,5,6,7\}$.  Given the colors used, and the fact that this is
to be a $9$-coloring, the corresponding coloring monomial is
$$
x_\emptyset^3\cd x_{36}^2\cd x_{367}\cd x_{[7]}^2.
$$

We record this in the following theorem.

\begin{theorem}%
There is a one-to-one correspondence between the $(n+1)$-colorings of $G$
and the monomials of degree $n$ in $K_G$.
\qed{} 
\end{theorem}
As an immediate consequence we now get the following.
\begin{corollary}\label{hilb-coro}%
  Let $F(K_G,t)$ be the Hilbert series of $K_G$ and let $d$ be the
  number of vertices in $G$.  Then
$$
F(K_G,t) = {1\over t}\cd\frac{W_G(t)}{(1-t)^{d+1}},
$$
Equivalently, $H(K_G,n)=\kig(n+1)$, that is, the Hilbert polynomial
of $K_G$ equals, up to a shift by one, the chromatic polynomial of
$G$.  \qed{}
\end{corollary}

\section{The coloring complex and its face ring}

Clearly, the basic coloring monomials are square-free.  Thus, the
quotient $R/K_G$ is the face ring of a simplicial complex whose vertex
set is a subset of $\{S\st S\sseq [d]\}$ and whose minimal non-faces
are 
$$
\{\; \{S_1,S_2,\ldots,S_k\} \st x_{S_1}x_{S_2}\cdots x_{S_k}
  \mbox{ is a basic coloring monomial}\}.
$$
We call this complex \emm the coloring complex of $G$, and denote
it by \dg.

One of the fundamental facts in the theory of face rings is that the
$h$-vector of a $d$-dimensional complex $\Delta$ is given by the
coefficients of the numerator when the Hilbert series of the face ring
$S$ of $\Delta$ is written as a rational function with denominator
$(1-t)^{d+1}$.  Thus, one customarily writes the Hilbert series
$F(S,t)$ of such a ring in the form
$$
F(S,t) = \frac{h(S,t)}{(1-t)^{d+1}} = 
\frac{h_0 +h_1t+\cdots+h_{d+1}t^{d+1}}{(1-t)^{d+1}},
$$
where $(h_0,h_1,\ldots,h_{d+1})$ is the $h$-vector of $\Delta$.

The following is easy to prove.
\begin{lemma}\label{diff-lemma}%
If $M$ is an ideal in a ring $S$ then
$F(S/M,t) = F(S,t) - F(M,t)$.
\qed
\end{lemma}

It is shown in Theorem \ref{struct-theorem} that the complex $\dg$,
whose face ring is $R/K_G$, has codimension one as a subcomplex of the
order complex of the Boolean algebra on $d$ atoms, whose face ring is
$R$.  Thus, letting $S=R$ and $M=K_G$ in Lemma \ref{diff-lemma}, we
obtain
\begin{eqnarray}\label{eq1}
  \frac{h(R/K_G,t)}{(1-t)^{d}} = F(R/K_G,t) = F(R,t) - F(K_G,t).
\end{eqnarray} 

Also, the ring $R$ is the face ring of a cone over the first
barycentric subdivision of a $(d-1)$-simplex and it is well-known that
its Hilbert series is given by
$$
F(R,t) = \frac{1}{t}\cd\frac{A_d(t)}{(1-t)^{d+1}},
$$
where $A_d(t)$ is the $d$-th Eulerian polynomial, which satisfies
$$
\sum_{n\ge0}{n^d t^n} = \frac{A_d(t)}{(1-t)^{d+1}}.
$$
Using this, together with identity \pref{eq1} and Corollary
\ref{hilb-coro}, allows us to relate the Hilbert series of $R/K_G$ to
the \emm tail,\/ of the chromatic polynomial of $G$, which we now
define.

\begin{definition}\rm%
  The \emm tail, of the chromatic polynomial $\ki_G$ of $G$ is
  $\tg(n)= n^d-\ki_G(n)$.
\end{definition} 

\begin{theorem}%
We have
$$
{1\over t}\sum_{n\ge 0}{\tg(n)t^n} = \frac{h(R/K_G,t)}{(1-t)^{d}}.
$$
Thus, up to a shift by one, the $W$-transform of the tail $\tg$ of
$\kig$ equals the polynomial whose coefficients are the coordinates of
the $h$-vector of the coloring complex of $G$ .
\end{theorem}
\pf{%
We have
\begin{eqnarray*}
{1\over t}\sum_{n\ge 0}{\tg(n)t^n} &=& 
\sum_{n\ge 0}{{(n+1)}^d t^n} - {1\over t}\sum_{n\ge 0}{\kig(n)t^n} =
\frac{1}{t}\cd\frac{A_d(t)}{(1-t)^{d+1}} - 
\frac{1}{t}\cd\frac{W_G(t)}{(1-t)^{d+1}}\\[6pt]
&=& F(R,t) - F(K_G,t) = F(R/K_G,t) = \frac{h(R/K,t)}{(1-t)^{d}}.
\end{eqnarray*} 

}\qed{} 

We shall now describe the structure of the complex \dg.  

First, however, to facilitate the following discussion we will let the
ring $R$ be the face ring of the order complex of the \emm truncated
Boolean algebra $\bdh$, on $d$ atoms, where $\bdh$ is $B_d$ with
$\emptyset$ and $[d]$ removed.  This is a harmless modification with
respect to our previous results (except in the rather trivial cases
when $G$ has fewer than three vertices) because the indeterminates
$x_{[d]}$ and $x_{\emptyset}$ divide none of the basic coloring
monomials.  Thus, the Hilbert series of $R/K_G$ is changed only in
that the denominator is divided by $(1-t)^2$.
The order complex of $\bdh$ is isomorphic to the first barycentric
subdivision of the boundary of a $(d-1)$-simplex, which has dimension
$d-2$.

This, of course, amounts to a redefinition of the coloring complex,
but even here the modification is trivial.  Namely, in the complex
\dg, $\emptyset$ and $[d]$ are both \emm cone points, that is, they
belong to every facet of \dg.
It is easy to show that removing a cone point (and all faces
containing it) from a complex changes the $h$-vector only by removing
its last coordinate, which is necessarily 0.  Thus, the $h$-vector of
$\dg$ remains essentially the same after removing $\emptyset$ and
$[d]$.

To each edge $e=ij$ of $G$, with $i<j$, we associate the $(d-1)!$
permutations of the letters in $S_e=[d]\sm\{i,j\}\cup \{e\}$. %
We call the permutations of $S_e$ $e$-permutations and we shall show
that the facets of $\dg$ correspond precisely to the edge permutations
for $G$.  It is important to note that, since we have removed the
cone point $[d]$ from our complex, the facet corresponding to an
edge-permutation does not contain the vertex $[d]$.  As an example,
the facet corresponding to the 25-permutation $3-25-1-4$ has vertices
$\{3\}$, $\{2,3,5\}$ and $\{1,2,3,5\}$.

\begin{theorem}\label{struct-theorem}%
  Let $G$ be a graph with $d\ge3$ vertices.
\begin{enumerate}%
\item%
  To each edge of $G$ there correspond exactly $(d-1)!$ facets of $\dg$ 
  and these are all the facets of \dg.  The sets of such facets for
  two distinct edges of $G$ are disjoint.  The facets thus
  corresponding to an edge form a $(d-3)$-sphere, which we call an
  \emm edge-sphere, and which is isomorphic to the order complex of a
  truncated Boolean algebra on $(d-1)$ elements.  That is, an
  edge-sphere is isomorphic to the first barycentric subdivision of
  the boundary of a $(d-3)$-simplex. In particular, $\dg$ has dimension
  $d-3$, unless $G$ is the graph with no edges, in which case $\dg$ is
  the empty complex.

\item\label{intersection}%
  Any two edge-spheres intersect in a ($d-4$)-sphere which is
  isomorphic to the order complex of a truncated Boolean algebra on
  $(d-2)$ elements.  Moreover, if $e$ and $f=ij$ are two edges, then
  the intersection of their two spheres separates the $e$-sphere into
  two halves, where one contains all vertices of $\dg$ that contain $i$
  and not $j$, whereas the other half contains those vertices that
  contain $j$ and not $i$.

\end{enumerate} 
\end{theorem}
\pf%
\begin{enumerate}%
\item\label{proof-1}%
 If $G$ has no edges, then the permutation $\pi= d\,(d-1)\ldots2\,1$
  has no cuts except 0, since $\ell(i)=0$ for all $i$.  Thus,
  $\pi$ corresponds to the empty monomial, or 1, and therefore the
  ideal $K_G$ is the entire ring $R$, so $\dg$ is the empty complex.  
  
  Suppose then that $e=ij$ is an edge in $G$.  Let $\pi = a_1
  a_2\cdots a_d$ where, for some $k$, we have $a_k=i$ and $a_{k+1}=j$.
  Let $S_0=\emptyset$ and let $F=\{S_1, S_2,\ldots,S_{d-2}\}$, where
$$
S_m=\cases{S_{m-1}\cup\{a_m\}, & if $m<k$,\cr
S_{m-1}\cup\{i,j\}, & if $m=k$,\cr
S_{m-1}\cup\{a_{m+1}\}, & if $m>k$.}
$$
Clearly, $F$ is a facet, since it has $d-2$ vertices and thus has
maximal dimension.  We claim that $\dg$ contains $F$.  If $\dg$ doesn't
contain $F$ then the ideal $K_G$ must contain a monomial dividing
{$x_{S_1}\cd x_{S_2}\cd \cdots\cd x_{S_{d-2}}$} and thus there must be
a permutation with a short $G$-sequence of which $S_1, S_2\sm S_1,
\ldots, S_{d-2}\sm S_{d-3}$ is a refinement.  But such a sequence must
contain a set containing both $i$ and $j$, which is a contradiction
since $ij$ is an edge in $G$ (see Lemma \ref{stable}).  We have thus
exhibited $(d-1)!$ facets $F$ associated to the edge $ij$. It is easy
to see that two facets thus corresponding to different edges are
different.  Namely, each vertex in such a facet is a set of vertices
from $G$ that either contains both or neither of a unique pair of
vertices in $G$, and this pair of vertices constitutes the edge
associated to the facet.
  
  Conversely, we need to show that any face of $\dg$ belongs to a facet
  associated to some edge of $G$ as described above.  Let
  $F=\{S_1\sseq S_2\sseq\cdots\sseq S_k\}$ be a face of \dg.  
  
  We first show that some of the difference sets $D_i=S_i\sm S_{i-1}$
  must contain the vertices of an edge in $G$.  If that is not the
  case then all the $D_i$ are stable sets in $G$. Construct a
  permutation $\pi$ in $\cls_d$ by first writing all the elements of
  $D_1$ in decreasing order, then those of $D_2$ in order of
  decreasing path lengths (w.r.t. to the vertices in $D_1$) and in
  decreasing order when two vertices have the same path length
  associated to them (see the proof of Theorem \ref{w-theorem}).
  Continue this way with all the $D_i$'s.  Then $D_1, D_2,\cdots, D_k$
  is a refinement of the short $G$-sequence of $\pi$, so the coloring
  ideal of $G$ contains a monomial dividing $x_{S_1} x_{S_2}\cdots
  x_{S_k}$.  This implies that $F$ is a non-face of \dg, a
  contradiction, so some of the $D_i$'s must contain an edge $e$ of
  $G$.
  
  This means that we can refine the chain of vertices of $F$ down to
  singletons except for having one of the sets in the refinement
  consist of the two vertices of $e$.  We have shown above that this
  refinement is a facet of \dg, and it is easy to see that it contains
  the face $F$.

\item%
  If the edges $e=ij$ and $f=km$ are disjoint, then the intersection
  of the $e$-sphere and the $f$-sphere consists of the subcomplex of
  $\dg$ whose vertices contain either both or neither of the vertices
  of the edge $e$ and, independently, either both or neither of the
  vertices of the edge $f$.  This subcomplex contains all faces of
  $\dg$ corresponding to permutations of $[d]$ where $i$ and $j$ are
  adjacent and in increasing order and where the same is true of $m$
  and $k$.
  
  If $e=ij$ and $f=kj$ are distinct edges then the intersection of the
  $e$-sphere and the $f$-sphere consists of the subcomplex of $\dg$ 
  whose vertices contain either all of $i,j,k$ or none of them.  This
  subcomplex contains all faces of $\dg$ corresponding to permutations
  of $[d]$ where $i$, $j$ and $k$ are three successive letters and in
  increasing order.
  
  In either case it is easy to see that the subcomplex of the
  intersection has dimension $d-4$ and that it separates the
  $e$-sphere as claimed.  Namely, there is a path along edges of \dg,
  not crossing the subcomplex, between any pair of vertices that
  belong to the same one of the halves described, but no such path
  between vertices in different halves.  \qed{}
\end{enumerate} 

\begin{remark}\rm\label{short-remark}%
  It is possible to consider the complex obtained in the same way as
  the coloring complex except that we don't remove the cone point
  $[d]$ from \dg.  This corresponds to associating the basic coloring
  monomials to the $G$-sequences of the permutations in $\cls_d$
  rather than their short $G$-sequences (and not stripping the face
  ring of the indeterminate $x_{[d]}$).  The $h$-vector of this
  complex is the $d$-th Eulerian vector (coefficients of the $d$-th
  Eulerian polynomial) plus the $h$-vector of the coloring complex
  shifted one step right. Thus, knowing the $h$-vector of this complex
  is equivalent to knowing the $h$ -vector of \dg.  As an example, for
  the graphs in Figure~\ref{fig2} we get $h$-vector
$$
(1,11,11,1) + (0,1,10,7) = (1,12,21,8),
$$
since the fourth Eulerian polynomial is $A_4(t)=t+11t^2+11t^3+t^4$.
\end{remark} 

Clearly, if the coloring complexes of two graphs are isomorphic, then
the graphs must be \emm chromatically equivalent,, that is, they must
have the same chromatic polynomial.  However, it is possible for two
non-isomorphic graphs to have isomorphic coloring complexes.  Namely,
there are two non-isomorphic graphs on $n\ge4$ vertices and two edges.
It follows from Theorem \ref{struct-theorem} that the coloring
complexes of these graphs must be isomorphic, because each consists of
two edge-spheres that intersect in a way independent of whether the
edges in question are disjoint.  Although these are the only examples
we know of non-isomorphic graphs with isomorphic coloring complexes we
do not know what the situation is in general.  In Figure \ref{fig2} we
give the coloring complexes of the two non-isomorphic --- but
chromatically equivalent --- graphs on three edges and four vertices.
As can be seen, these complexes are not isomorphic (one has a
``triangle'' and the other one doesn't).

\setlength{\unitlength}{.8mm}
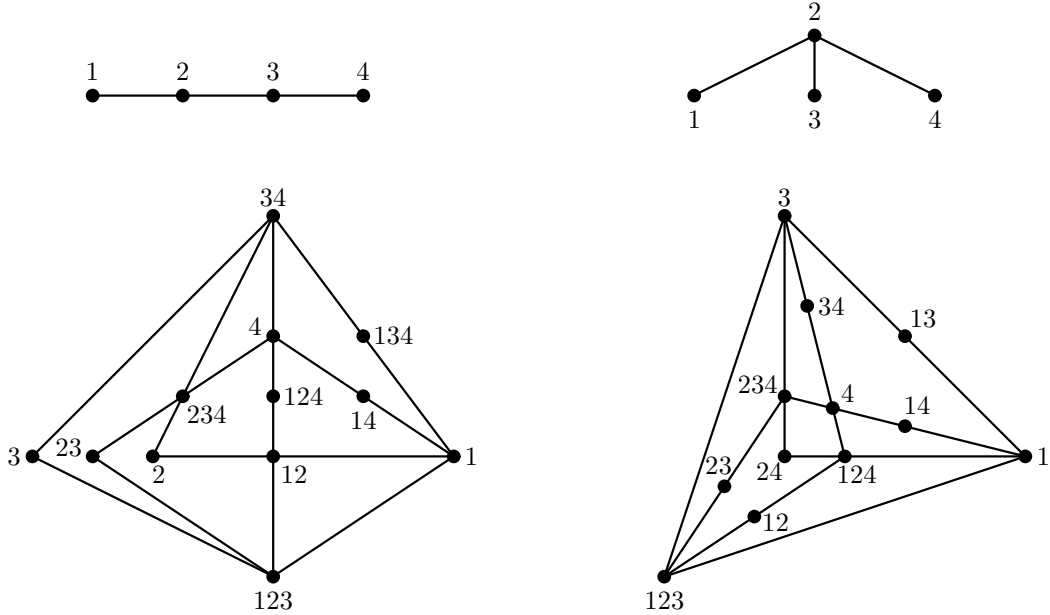
\begin{figure}[t]
\begin{center}
\begin{picture}(200,117)(5,-40)
\newcommand\pt{\circle*{2}}

\thicklines

\small

\put(20,60){
\put(0,0){\pt\makebox(0,8){1}}
\put(15,0){\pt\makebox(0,8){2}}
\put(30,0){\pt\makebox(0,8){3}}
\put(45,0){\pt\makebox(0,8){4}}

\path(0,0)(15,0)

\path(15,0)(30,0)

\path(30,0)(45,0)
}

\put(50,0){\pt\makebox(7,-6){12}}
\put(50,-20){\pt\makebox(0,-8){123}}
\put(10,0){\pt\makebox(-6,0){3}}
\put(50,40){\pt\makebox(0,6){34}}
\put(50,20){\pt\makebox(-6,3){4}}
\put(50,10){\pt\makebox(10,0){124}}

\put(20,0){\pt\makebox(-8,3){23}}
\put(35,10){\pt\makebox(8,-6){234}}
\put(65,10){\pt\makebox(0,-8){14}}
\put(80,0){\pt\makebox(6,0){1}}

\put(30,0){\pt\makebox(2,-6){2}}
\put(65,20){\pt\makebox(10,0){134}}

\path(80,0)(50,0)(30,0)(35,10)(50,40)(65,20)(80,0)

\path(50,-20)(20,0)(35,10)(50,20)(65,10)(80,0)(50,-20)

\path(50,-20)(10,0)(50,40)(50,20)(50,10)(50,0)(50,-20)

\put(115,-20){%

\put(5,90){
\put(20,0){\pt\makebox(0,8){2}}
\put(0,-10){\pt\makebox(0,-8){1}}
\put(20,-10){\pt\makebox(0,-8){3}}
\put(40,-10){\pt\makebox(0,-8){4}}

\path(0,-10)(20,0)

\path(20,0)(20,-10)

\path(20,0)(40,-10)
}

\put(20,20){\pt\makebox(-5,-5){24}}
\put(20,30){\pt\makebox(-9,4){234}}
\put(20,60){\pt\makebox(0,6){3}}
\put(40,40){\pt\makebox(6,6){13}}
\put(60,20){\pt\makebox(6,0){1}}
\put(30,20){\pt\makebox(4,-6){124}}

\path(20,20)(20,60)(60,20)(20,20)

\put(0,0){\pt\makebox(0,-8){123}}
\put(23.75,45){\pt\makebox(8,0){34}}
\put(28,28){\pt\makebox(5,5){4}}
\put(15,10){\pt\makebox(7,-2){12}}

\path(0,0)(20,60)(30,20)(0,0)

\put(10,15){\pt\makebox(-2,7){23}}
\put(40,25){\pt\makebox(4,7){14}}

\path(0,0)(20,30)(60,20)(0,0)

}

\end{picture}
\caption{\label{fig2} Two non-isomorphic graphs (trees) with the same chromatic
  polynomial but non-isomorphic coloring complexes.
}
\end{center}
\end{figure}

Perhaps more interesting is that the coloring complex distinguishes
the two graphs given in Stanley's paper \cite[Figure 1]{stan-symm},
which his chromatic symmetric function does not distinguish (see
Figure \ref{fig3}).  This can be seen as follows: A complex on a given
vertex set is determined by (in fact equivalent to) its set of minimal
non-faces, which in turn is equivalent to the (unique) minimal set of
generators of the ideal defining its face ring.  Suppose the graphs
$G$ and $H$ in Figure \ref{fig3} have isomorphic coloring complexes.
Then there is a bijection $\phi$ between their vertex sets so that $A$
is a minimal non-face of $\dg$ if and only if $\phi(A)$ is a minimal
non-face of $\Delta_H$.  Thus, the minimal sets of generators for the
coloring ideals of $G$ and $H$ must have the same number of monomials
of each degree.  Also, the multiplicities of corresponding
indeterminates in the sets of monomials constituting the respective
minimal generating sets for each coloring ideal must be the same.
This is not the case for the coloring ideals of $G$ and $H$, which we
have verified with the aid of the computer algebra program {\sc
  Macaulay 2} \cite{mac2}.

\setlength{\unitlength}{1mm}

\begin{figure}
\begin{picture}(150,35)(-13,-10)
\newcommand\pt{\circle*{1.5}}

\path(0,0)(0,20)(40,0)(40,20)(0,0)

\put(0,0){\pt}
\put(0,20){\pt}
\put(20,10){\pt}
\put(40,0){\pt}
\put(40,20){\pt}

\put(80,0){
\path(0,20)(20,20)(20,0)(0,0)(0,20)(20,0)(20,20)(30,0)

\put(0,0){\pt}
\put(0,20){\pt}
\put(20,20){\pt}
\put(20,0){\pt}
\put(30,0){\pt}
}

\end{picture}
\caption{\label{fig3} Another two non-isomorphic graphs with the same chromatic
  polynomial and same symmetric chromatic function but non-isomorphic
  coloring complexes. (See \protect{\cite[Figure 1]{stan-symm}}.)}
\end{figure}
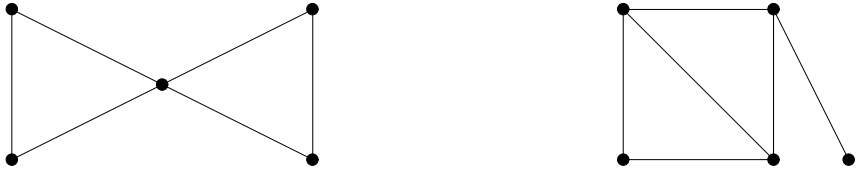

The next corollary follows from part \ref{proof-1} of Theorem
\ref{struct-theorem}.  It can also be proved directly from the
well-known fact that the coefficient to $-n^{d-1}$ in \kig, which is
the leading coefficient of \tg, equals the number of edges in $G$.

\begin{corollary}\label{numfacets}%
  The number of facets of \dg, and thus the sum of the coefficients of
  the $h$-vector of \dg, is $E\cd(d-1)!$, where $E$ is the number of
  edges in $G$.  
\qed{}
\end{corollary}

\begin{theorem}%
  The Euler characteristic of  $\dg$ equals the number of acyclic
  orientations of $G$.
\end{theorem}
\pf{%
  Up to a sign, the \emm reduced,\/ Euler characteristic of a
  $(d-1)$-dimensional complex $\Delta$ is equal to the $d$-th
  coordinate $h_{d}(\Delta)$ of the $h$-vector of $\Delta$ and thus
  the reduced Euler characteristic of $\dg$ equals the leading
  coefficient of the $W$-transform of the tail $\tg(n)$ of the
  chromatic polynomial.  It is well known (and easy to prove) that the
  leading coefficient of the $W$-transform of a polynomial $P(x) = a_0
  +a_1 x +\cdots +a_d x^d$ equals, up to a sign, the alternating sum
  of the coefficients of $P$.  More precisely, it equals $(-1)^d
  P(-1)$, where $d$ is the degree of $P$.  Clearly, $(-1)^{d-1} T(-1)
  = (-1)^d\kig(-1) -1$.  But, by a theorem of Stanley \cite[Corollary
  1.3]{stan-chrom}, $(-1)^d\kig(-1)$ equals the number of acyclic
  orientations of $G$.  Since the reduced Euler characteristic is one
  less than the Euler characteristic, this establishes the claim.
\qed{}

\centerline{\sc Open problems}
\bigskip

An obvious question is whether $\dg$ is shellable.  A consequence of
shellability would be that $\dg$ (equivalently $R$) is Cohen-Macaulay
In that case the $h$-vector of the coloring complex must be an \emm
$M$-vector, (see \cite{CCA}), which would put certain restrictions on
the values of the tail $T_G$ and thereby on the values of the
chromatic polynomial.

For $i=1,\ldots,d$, let $\theta_i=\sum_{|S|=i}{x_S}$.  Then it can be
shown that $\theta=\theta_1,\theta_2,\ldots,\theta_d$ is a homogeneous
(linear) system of parameters for $R$.  Is $R$ a free
$k[\theta]$-module?  That is equivalent to $R$ being Cohen-Macaulay.
However, a proof of this would likely be equivalent to finding a
shelling of \dg, and shellability of $\dg$ would imply that $\dg$
(equivalently $R$) is Cohen-Macaulay.

It might be interesting to know what the minimal set of generators is
for the coloring ideal of a graph $G$ and in particular what the size
of this set is.  Perhaps it is more interesting to determine this for
the ideal $K\cup I$, where $K$ is the coloring ideal of $G$ and $I$ is
the ideal used in defining the ring $R$, since the face ring of the
coloring complex $C$ is given by $A/(K\cup I)$.

\bigskip\bigskip

\centerline{\sc Acknowledgements}
\smallskip

My thanks to Eric Babson, to Ralf Fröberg, to Richard Stanley, and to
an anonymous refereee for valuable comments.

\end{document}